\documentclass[a4paper,12pt]{article}

\usepackage{amsthm}
\usepackage{amsmath,amssymb,latexsym,amsfonts,mathrsfs}

\theoremstyle{plain}
\newtheorem{Thm}{Theorem}[section]

\newtheorem{Lem}[Thm]{Lemma}
\newtheorem{Prop}[Thm]{Proposition}

\theoremstyle{definition}
\newtheorem{Def}[Thm]{Definition}

\newcommand{\Proof}[2][Proof]{\begin{proof}[{#1}] #2 \end{proof}}

\numberwithin{equation}{section}

\makeatletter
\renewcommand\section{\@startsection {section}{1}{\z@}%
                                   {-3.5ex \@plus -1ex \@minus -.2ex}%
                                   {2.3ex \@plus.2ex}%
                                   {\normalfont\large\bf}}
\makeatother

\makeatletter
\renewcommand\subsection{\@startsection {subsection}{1}{\z@}%
                                   {-3.5ex \@plus -1ex \@minus -.2ex}%
                                   {2.3ex \@plus.2ex}%
                                   {\normalfont\normalsize\bf}}
\makeatother

\makeatletter
\@addtoreset{footnote}{page}
\makeatother

\newcommand{\rbra}[1]{\left( #1 \right)} 
\newcommand{\cbra}[1]{\left\{ #1 \right\}} 
\newcommand{\sbra}[1]{\left[ #1 \right]} 
\newcommand{\abra}[1]{\left\langle #1 \right\rangle} 

\newcommand{\bmat}[1]{\begin{bmatrix} #1 \end{bmatrix}} 

\renewcommand{\subitem}{\par\hangindent 5mm \hspace*{-2mm}}
\renewcommand{\subsubitem}{\par\hangindent 10mm \hspace*{2mm}}

\renewcommand{\d}{{\rm d}} 
\newcommand{\dist}{\stackrel{{\rm d}}{=}}
\newcommand{\tend}[2]{\mathrel{\mathop{\longrightarrow}\limits^{#1}_{#2}}}
\renewcommand{\hat}{\widehat}
\renewcommand{\tilde}{\widetilde}
\newcommand{\Supp}{\mathop{\rm Supp}}

\newcommand{\bN}{\ensuremath{\mathbb{N}}}
\newcommand{\bR}{\ensuremath{\mathbb{R}}}
\newcommand{\bZ}{\ensuremath{\mathbb{Z}}}

\newcommand{\cF}{\ensuremath{\mathcal{F}}}
\newcommand{\fS}{\ensuremath{\mathfrak{S}}}

\newcommand{\minip}[2]{
\begin{minipage}{#1}
\begin{center}
#2
\end{center}
\end{minipage}}

\begin{document}
\begin{center}
{\Large \bf Random walk in a finite directed graph subject to a road coloring}
\end{center}

\begin{center}
Kouji Yano\footnote{
Graduate School of Science, Kobe University, Kobe, JAPAN.}\footnote{
The research of this author was supported by KAKENHI (20740060).} 
\end{center}

\begin{center}
{\small \today}
\end{center}

\begin{abstract}
A necessary and sufficient condition for 
a random walk in a finite directed graph subject to a road coloring 
to be measurable with respect to the driving random road colors 
is proved to be that the road coloring is synchronizing. 
For this, the random walk subject to a non-synchronizing road coloring is proved 
to have uniform law on a certain partition of the state space. 
\end{abstract}


\section{Introduction}

Let us consider a finite directed graph of constant outdegree. 
See, for example, Figure 1 below; 
there are five sites 
and from each site there are two oneway roads laid. 
Let us color every road blue or red 
so that no two roads running from the same site have the same color. 
See, for example, Figure 2 below; 
the thick roads are colored red and the thin ones blue. 

\minip{7cm}{
\unitlength 0.1in
\begin{picture}( 20.0000, 18.0000)(  0.0000,-18.0000)
%
\special{pn 8}%
\special{ar 200 800 200 200  0.0000000 6.2831853}%
%
\special{pn 8}%
\special{ar 1000 200 200 200  0.0000000 6.2831853}%
%
\special{pn 8}%
\special{ar 1800 800 200 200  0.0000000 6.2831853}%
%
\special{pn 8}%
\special{ar 600 1600 200 200  0.0000000 6.2831853}%
%
\special{pn 8}%
\special{ar 1400 1600 200 200  0.0000000 6.2831853}%
\put(10.0000,-2.0000){\makebox(0,0){1}}%
\put(2.0000,-8.0000){\makebox(0,0){2}}%
\put(6.0000,-16.0000){\makebox(0,0){3}}%
\put(14.0000,-16.0000){\makebox(0,0){4}}%
\put(18.0000,-8.0000){\makebox(0,0){5}}%
%
\special{pn 8}
\special{ar 1800 550 184 184  2.4858970 6.2831853}%
\special{ar 1800 550 184 184  0.0000000 0.6202495}%
\special{pa 1550 610}%
\special{pa 1650 660}%
\special{fp}%
\special{pa 1640 640}%
\special{pa 1700 550}%
\special{fp}

\special{pn 8}%
\special{pa 790 230}%
\special{pa 280 610}%
\special{fp}%
\special{pa 280 500}%
\special{pa 280 610}%
\special{fp}%
\special{pa 280 610}%
\special{pa 400 610}%
\special{fp}%
%
\special{pn 8}%
\special{pa 390 810}%
\special{pa 1590 810}%
\special{fp}%
\special{pa 1470 730}%
\special{pa 1590 800}%
\special{fp}%
\special{pa 1590 800}%
\special{pa 1460 900}%
\special{fp}%
%
\special{pn 8}%
\special{pa 750 1450}%
\special{pa 1630 920}%
\special{fp}%
\special{pa 1480 930}%
\special{pa 1640 920}%
\special{fp}%
\special{pa 1640 920}%
\special{pa 1550 1060}%
\special{fp}%
%
\special{pn 8}%
\special{pa 1360 1410}%
\special{pa 1070 380}%
\special{fp}%
\special{pa 1030 500}%
\special{pa 1060 400}%
\special{fp}%
\special{pa 1060 400}%
\special{pa 1190 460}%
\special{fp}%
%
\special{pn 8}%
\special{pa 1720 990}%
\special{pa 1500 1420}%
\special{fp}%
\special{pa 1470 1280}%
\special{pa 1490 1430}%
\special{fp}%
\special{pa 1490 1430}%
\special{pa 1640 1380}%
\special{fp}%
%
\special{pn 8}%
\special{pa 850 330}%
\special{pa 380 710}%
\special{fp}%
\special{pa 390 580}%
\special{pa 390 690}%
\special{fp}%
\special{pa 390 690}%
\special{pa 520 710}%
\special{fp}%
%
\special{pn 8}%
\special{pa 270 980}%
\special{pa 500 1430}%
\special{fp}%
\special{pa 360 1370}%
\special{pa 490 1410}%
\special{fp}%
\special{pa 490 1410}%
\special{pa 520 1260}%
\special{fp}%
%
\special{pn 8}%
\special{pa 790 1610}%
\special{pa 1200 1610}%
\special{fp}%
\special{pa 1110 1530}%
\special{pa 1200 1610}%
\special{fp}%
\special{pa 1200 1610}%
\special{pa 1110 1700}%
\special{fp}%
%
\special{pn 8}%
\special{pa 1250 1470}%
\special{pa 360 910}%
\special{fp}%
\special{pa 400 1060}%
\special{pa 380 920}%
\special{fp}%
\special{pa 380 920}%
\special{pa 550 910}%
\special{fp}%
\end{picture}%
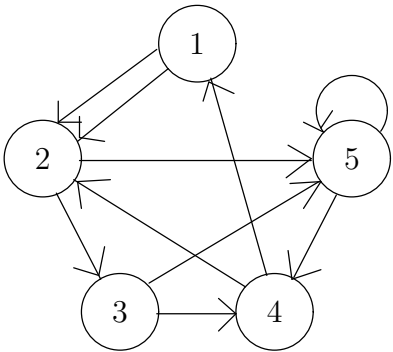
\\ Figure 1. 
}
\minip{7cm}{
\unitlength 0.1in
\begin{picture}( 20.0000, 18.0000)(  0.0000,-18.0000)
%
\special{pn 8}%
\special{ar 200 800 200 200  0.0000000 6.2831853}%
%
\special{pn 8}%
\special{ar 1000 200 200 200  0.0000000 6.2831853}%
%
\special{pn 8}%
\special{ar 1800 800 200 200  0.0000000 6.2831853}%
%
\special{pn 8}%
\special{ar 600 1600 200 200  0.0000000 6.2831853}%
%
\special{pn 8}%
\special{ar 1400 1600 200 200  0.0000000 6.2831853}%
\put(10.0000,-2.0000){\makebox(0,0){1}}%
\put(2.0000,-8.0000){\makebox(0,0){2}}%
\put(6.0000,-16.0000){\makebox(0,0){3}}%
\put(14.0000,-16.0000){\makebox(0,0){4}}%
\put(18.0000,-8.0000){\makebox(0,0){5}}%
%
\special{pn 20}
\special{ar 1800 550 184 184  2.4858970 6.2831853}%
\special{ar 1800 550 184 184  0.0000000 0.6202495}%
\special{pa 1550 610}%
\special{pa 1650 660}%
\special{fp}%
\special{pa 1640 640}%
\special{pa 1700 550}%
\special{fp}

\special{pn 8}%
\special{pa 790 230}%
\special{pa 280 610}%
\special{fp}%
\special{pa 280 500}%
\special{pa 280 610}%
\special{fp}%
\special{pa 280 610}%
\special{pa 400 610}%
\special{fp}%
%
\special{pn 8}%
\special{pa 390 810}%
\special{pa 1590 810}%
\special{fp}%
\special{pa 1470 730}%
\special{pa 1590 800}%
\special{fp}%
\special{pa 1590 800}%
\special{pa 1460 900}%
\special{fp}%
%
\special{pn 8}%
\special{pa 750 1450}%
\special{pa 1630 920}%
\special{fp}%
\special{pa 1480 930}%
\special{pa 1640 920}%
\special{fp}%
\special{pa 1640 920}%
\special{pa 1550 1060}%
\special{fp}%
%
\special{pn 8}%
\special{pa 1360 1410}%
\special{pa 1070 380}%
\special{fp}%
\special{pa 1030 500}%
\special{pa 1060 400}%
\special{fp}%
\special{pa 1060 400}%
\special{pa 1190 460}%
\special{fp}%
%
\special{pn 8}%
\special{pa 1720 990}%
\special{pa 1500 1420}%
\special{fp}%
\special{pa 1470 1280}%
\special{pa 1490 1430}%
\special{fp}%
\special{pa 1490 1430}%
\special{pa 1640 1380}%
\special{fp}%
%
\special{pn 20}%
\special{pa 850 330}%
\special{pa 380 710}%
\special{fp}%
\special{pa 390 580}%
\special{pa 390 690}%
\special{fp}%
\special{pa 390 690}%
\special{pa 520 710}%
\special{fp}%
%
\special{pn 20}%
\special{pa 270 980}%
\special{pa 500 1430}%
\special{fp}%
\special{pa 360 1370}%
\special{pa 490 1410}%
\special{fp}%
\special{pa 490 1410}%
\special{pa 520 1260}%
\special{fp}%
%
\special{pn 20}%
\special{pa 790 1610}%
\special{pa 1200 1610}%
\special{fp}%
\special{pa 1110 1530}%
\special{pa 1200 1610}%
\special{fp}%
\special{pa 1200 1610}%
\special{pa 1110 1700}%
\special{fp}%
%
\special{pn 20}%
\special{pa 1250 1470}%
\special{pa 360 910}%
\special{fp}%
\special{pa 400 1060}%
\special{pa 380 920}%
\special{fp}%
\special{pa 380 920}%
\special{pa 550 910}%
\special{fp}%
\end{picture}%
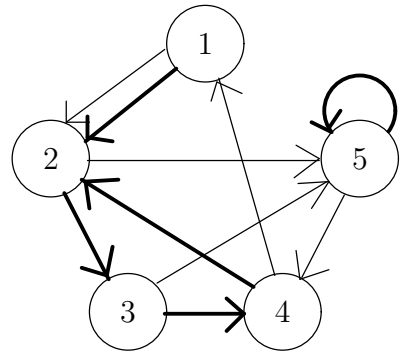
\\ Figure 2. 
}

We call $ N $ a {\em random color} 
if $ N $ is a random variable which takes values in the set of road colors; 
in this case, red and blue. 
To put it roughly, we mean by {\em random walk} 
a pair of processes $ \{ X,N \} $ 
where $ N=(N_k)_{k \in \bZ} $ is a sequence of colors 
which are independent and identically distributed (abbreviated as IID) 
and $ X=(X_k)_{k \in \bZ} $ is a site-valued process 
which moves at each step from $ X_{k-1} $ to $ X_k $ 
being driven by the random color $ N_k $. 
We have the following table for instance: 

\begin{center}
\begin{tabular}{r||c|c|c|c|c|c}
$ k $   & $ k_0-5 $ & $ k_0-4 $ & $ k_0-3 $ & $ k_0-2 $ & $ k_0-1 $ & $ k_0 $ \\
\hline
\hline
$ N_k $ &    & blue & red & blue & blue & red \\
\hline
$ X_k $ & 1  & 2    & 3   & 5    & 4    & 2 \\
\hline
$ X_k $ & 2  & 5    & 5   & 4    & 1    & 2 \\
\hline
$ X_k $ & 3  & 5    & 5   & 4    & 1    & 2 \\
\hline
$ X_k $ & 4  & 1    & 2   & 5    & 4    & 2 \\
\hline
$ X_k $ & 5  & 4    & 2   & 5    & 4    & 2 
\end{tabular}
\end{center}

In this table, we let $ k_0 $ be a certain time 
and we assume that 
at the 5 steps to the time $ k=k_0 $ 
the random colors are blue, red, blue, blue and red, in this order. 
Then, no matter how the process $ X $ moves before $ k_0-5 $, 
the value of $ X_{k_0} $ falls 2, 
and therefore the values of $ X $ afterward can be known from the values of $ N $. 
We are interested in the necessary and sufficient condition 
that we can always know the values of $ X $ only from the values of $ N $.

\subsection{Random walk subject to a road coloring}

Let $ V $ be a set of finite symbols. 
A matrix $ [A(y,x)]_{y,x \in V} $ whose entries are non-negative integers 
is called an {\em adjacency matrix} 
and the pair $ (V,A) $ a {\em directed graph}. 
Note that there may be multiple edges 
which are not distinguished from each other. 
We call each element of $ V $ a {\em site}. 
For each $ y,x \in V $, the value of $ A(y,x) $ 
may be regarded as the number of {\em (oneway) roads} from $ x $ to $ y $. 
(We prefer to write $ A(y,x) $ than write $ A(x,y) $.) 

From each site $ x \in V $ there are as many roads as 
$ \deg(x) := \sum_{y \in V} A(y,x) $. 
The graph $ (V,A) $ or the adjacency matrix $ A $ is called {\em $ d $-out} 
if $ \deg(x) \equiv d $; in other words, there are $ d $ roads from each site. 
It is called {\em of constant outdegree} if it is $ d $-out for some $ d $. 
Write 
\begin{align}
\Sigma = \text{the set of all mappings from $ V $ to itself}. 
\label{}
\end{align}
For $ \sigma_1,\sigma_2 \in \Sigma $ and $ x \in V $, 
we write $ \sigma_2 \sigma_1 x $ simply for $ \sigma_2 ( \sigma_1 (x) ) $. 
The semigroup $ \Sigma $ acts on $ V $ as follows: 
\begin{align}
(\sigma_1 \sigma_2) x 
= \sigma_1 (\sigma_2 x) 
, \quad \sigma_1,\sigma_2 \in \Sigma , \ x \in V . 
\label{eq: action}
\end{align}
Each element $ \sigma \in \Sigma $ may be identified with 
the 1-out adjacency matrix $ [\sigma(y,x)]_{y,x \in V} $ via equation 
\begin{align}
\sigma(y,x) = 1_{\{ y=\sigma x \}} . 
\label{}
\end{align}
An adjacency matrix $ A $ of constant outdegree 
admits a family $ C = \{ \sigma^{(1)},\ldots,\sigma^{(d)} \} $ of $ \Sigma $ 
(possibly with repeated elements) such that 
\begin{align}
A = \sigma^{(1)} + \cdots + \sigma^{(d)} . 
\label{eq: adjacency}
\end{align}
Such a family $ C $ will be called a {\em road coloring}, 
because $ C $ indicates one of the ways 
of coloring the $ d $ roads running from each site 
so that no two roads from the same site have the same color. 

Let $ \mu $ be a probability law on $ \Sigma $. 
We write $ \Supp(\mu) $ for the support of $ \mu $, i.e., 
\begin{align}
\Supp(\mu) = \cbra{ \sigma \in \Sigma : \mu(\sigma)>0 } . 
\label{}
\end{align}
Enumerating $ \Supp(\mu) $ as $ \{ \sigma^{(1)},\ldots,\sigma^{(d)} \} $, 
we define the adjacency matrix $ A $ by \eqref{eq: adjacency} 
so that $ \Supp(\mu) $ is a road coloring of $ (V,A) $. 
The resulting directed graph $ (V,A) $ is called 
the {\em directed graph induced by $ \mu $}. 
Now we introduce random walk in a directed graph 
indexed by $ \bZ := \{ \ldots,-1,0,1,\ldots \} $ as follows. 

\begin{Def}
Let $ \mu $ be a probability law on $ \Sigma $. 
A pair of processes $ \{ X,N \} $ defined on a probability space 
is called a {\em $ \mu $-random walk} 
if $ X=(X_k)_{k \in \bZ} $ and $ N=(N_k)_{k \in \bZ} $ are processes 
taking values in $ V $ and $ \Sigma $, respectively, 
such that the following statements hold: 
\begin{enumerate}
\item 
$ N_k $ is independent of $ \sigma(X_j,N_j:j \le k-1) $ for each $ k \in \bZ $; 
\item 
$ N=(N_k)_{k \in \bZ} $ is IID with common law $ \mu $; 
\item it holds that 
\begin{align}
X_k = N_k X_{k-1} 
, \quad \text{a.s.}, \ k \in \bZ . 
\label{eq: evolution}
\end{align}
\end{enumerate}
\end{Def}

Let $ (V,A) $ denote the directed graph induced by $ \mu $. 
The process $ X=(X_k)_{k \in \bZ} $ moves at each step from a site to another in $ (V,A) $, 
being driven by the randomly-chosen road colors indicated by $ N=(N_k)_{k \in \bZ} $ 
via equation \eqref{eq: evolution}. 
This is why we call such a process $ \{ X,N \} $ a random walk. 
Note that our definition is different from the one in a lot of literatures; 
see \cite{MR1743100} and references therein. 

\begin{Def}
A $ \mu $-random walk $ \{ X,N \} $ is called {\em strong} if 
$ X_k $ is a.s. measurable with respect to $ \sigma(N_j:j \le k) $ for all $ k \in \bZ $; 
or equivalently, there exist measurable mappings 
$ f_k:\Sigma^{\bN} \to V $ such that 
\begin{align}
X_k = f_k(N_k,N_{k-1},\ldots) 
\quad \text{a.s. for all $ k \in \bZ $}. 
\label{}
\end{align}
\end{Def}

The purpose of this paper is to investigate a necessary and sufficient condition 
for the $ \mu $-random walk to be strong.

\subsection{Main theorem}

Let $ A^0 $ denote the identity matrix 
and define $ A^n(y,x) = \sum_{z \in V} A^{n-1}(y,z) A(z,x) $, $ y,x \in V $ 
recursively for $ n \ge 1 $. 
A directed graph $ (V,A) $ is called {\em strongly-connected} 
if for any $ y,x \in V $ there exists $ n=n(y,x) \ge 1 $ such that $ A^n(y,x) \ge 1 $; 
or in other words, one can walk from every site to every other site. 
The graph $ (V,A) $ is called {\em aperiodic} 
if the {\em period} at $ x \in V $, 
i.e., the greatest common divisor among $ \{ n \ge 1 : A^n(x,x) \ge 1 \} $, 
is one for all $ x \in V $. 
Note that $ (V,A) $ is both strongly connected and aperiodic if and only if 
there exists a positive integer $ r $ such that 
$ A^r(y,x) \ge 1 $ for all $ y,x \in V $. 
We say that the directed graph $ (V,A) $ {\em satisfies the assumption} {\bf (A)} 
if it is of constant outdegree, strongly-connected, and aperiodic. 
We will prove as Theorem \ref{thm: main1} that 
under the assumption {\bf (A)} 
there exists a unique $ \mu $-random walk $ \{ X,N \} $ which is stationary. 

For $ \sigma \in \Sigma $, we write $ \sigma V = \{ \sigma x: x \in V \} $. 
A road coloring $ C $ is called {\em synchronizing} if 
there exists a sequence $ s=(\sigma_p,\ldots,\sigma_1) $ of road colors 
such that the composition $ \abra{s} := \sigma_p \cdots \sigma_1 $ 
maps $ V $ onto a singleton; 
or in other words, 
those who walk in the directed graph being driven by the road colors 
$ \sigma_1,\ldots,\sigma_p $ in this order will arrive at a common site. 

Now one of our main results is as follows. 

\begin{Thm} \label{thm: main intro}
Let $ \mu $ be a probability law on $ \Sigma $ 
and let $ \{ X,N \} $ be a $ \mu $-random walk. 
Suppose that the directed graph induced by $ \mu $ satisfies the assumption {\bf (A)}. 
Then the following three assertions are equivalent: 
\subitem {\rm (i)} 
$ \Supp(\mu) $ is synchronizing. 
\subitem {\rm (ii)} 
The limit $ \displaystyle \lim_{l \to -\infty } N_k N_{k-1} \cdots N_{l+1} $ 
exists a.s. for all $ k \in \bZ $. 
\subitem {\rm (iii)} 
The $ \mu $-random walk $ \{ X,N \} $ is strong.  
\end{Thm}

Theorem \ref{thm: main intro} will be proved in Section \ref{sec: main2}. 
Note that the most difficult part of Theorem \ref{thm: main intro} 
is to show that (iii) implies (i). 
We shall prove the contraposition: 
\begin{align}
\text{If $ \Supp(\mu) $ is non-synchronizing, 
then the $ \mu $-random walk $ \{ X,N \} $ is non-strong.} 
\label{eq: main2}
\end{align}

\subsection{A typical sufficient condition}

To prove \eqref{eq: main2}, we need to find some extra randomness 
which is not measurable with respect to $ \sigma(N_j:j \le k) $. 
The key to the proof is to reveal a certain symmetry, 
or to put it more precisely, 
to construct another random walk from the original random walk 
by stopping it at certain stopping times 
and the problem is then reduced to the proposition given as follows. 

Let $ \fS(V) $ denote the permutation group of $ V $, 
which may be regarded as a subgroup of $ \Sigma $. 

\begin{Prop} \label{prop: typical}
Let $ \{ X,N \} $ be a $ \mu $-random walk. 
Suppose that $ \Supp(\mu) $ is contained in $ \fS(V) $ and that 
\begin{align}
\mu^{*n} \rbra{ \sigma \in \Sigma : \sigma(i)=j } 
\tend{}{n \to \infty } 
\frac{1}{\sharp(V)} 
, \quad i,j \in V 
\label{}
\end{align}
where $ \sharp(V) $ stands for the number of elements of $ V $ 
and $ \mu^{*n} $ for the $ n $-times convolution of $ \mu $. 
Then 
the stationary law of $ X $ is uniform law on $ V $ 
and that $ \{ X,N \} $ is non-strong. 
\end{Prop}

\Proof{
Let $ \lambda $ denote the stationary law for the process $ X=(X_k)_{k \in \bZ} $. 
Since $ X_0 = N_0 N_{-1} \cdots N_{-n+1} X_{-n} $, we have 
\begin{align}
\lambda(j) 
= P(X_0=j) 
=& \sum_{i \in V} P(N_0 N_{-1} \cdots N_{-n+1}(i)=j) P(X_{-n}=i) 
\label{} \\
=& \sum_{i \in V} \mu^{*n} \rbra{ \sigma \in \Sigma : \sigma(i)=j } \lambda(i) 
\label{} \\
\tend{}{n \to \infty }& 
\frac{1}{\sharp(V)} \sum_{i \in V} \lambda(i) 
= \frac{1}{\sharp(V)} . 
\label{}
\end{align}
Let $ f $ be an arbitrary function on $ V $ 
and let $ l $ be a negative integer. 
Since $ X_0 = N_0 N_{-1} \cdots N_{l+1} X_l $, we have 
\begin{align}
E[f(X_0)|\sigma(N_j:j \ge l+1)] 
=& \left. E[f(\sigma X_l)] \right|_{\sigma = N_0 N_{-1} \cdots N_{l+1}} 
\label{} \\
=& \left. \int_V f(\sigma x) \lambda(\d x) \right|_{\sigma = N_0 N_{-1} \cdots N_{l+1}} 
\label{}
\end{align}
Since $ \Supp(\mu) \subset \fS(V) $ and since the uniform law on $ V $ is $ \fS(V) $-invariant, 
we see that 
\begin{align}
E[f(X_0)|\sigma(N_j:j \ge l+1)] 
= \int_V f(x) \lambda(\d x) = E[f(X_0)] . 
\label{}
\end{align}
Letting $ l \to -\infty $, we have 
\begin{align}
E[f(X_0)|\sigma(N_j:j \in \bZ)] = E[f(X_0)] . 
\label{}
\end{align}
This shows that $ X_0 $ is independent of $ \sigma(N_j:j \in \bZ) $. 
}

\subsection{Backgrounds}

Let us give a brief remark on the backgrounds of this study. 

{\bf a). Road coloring problem.} 
What we call the {\em road coloring problem} is the following: 
\begin{align}
\text{Does a directed graph admit at least one synchronizing road coloring?} 
\label{}
\end{align}
This problem was first posed 
by Adler--Goodwyn--Weiss \cite{MR0437715} (see also \cite{MR0257315}) 
in the context of the isomorphism problem of symbolic dynamics 
with a common topological entropy. 
It was solved only recently by Trahtman \cite{MR2534238}, 
so that the problem is now a theorem: 

\begin{Thm}[{Trahtman \cite{MR2534238}}]
A directed graph which satisfies the assumption {\bf (A)} 
admits at least one synchronizing road coloring. 
\end{Thm}

See also \cite{MR1788119}, \cite{MR2312963} and \cite{MR2081291} 
for some other developments before Trahtman \cite{MR2534238}. 


{\bf b). Tsirelson's equation in discrete time.} 
Stochastic equation \eqref{eq: evolution} is related to 
the study \cite{MR0375461} by Tsirelson, who introduced, 
in order to construct an example 
of a stochastic differential equation which has a non-strong solution, 
a stochastic equation indexed by the negative integer 
\begin{align}
X_k = N_k + X_{k-1} 
, \quad k \in -\bN , 
\label{eq: Tsirel eq}
\end{align}
where $ X $ and $ N $ take values in the one-dimensional torus $ \bR/\bZ $. 

Yor \cite{MR1147613} obtained a necessary and sufficient condition 
in terms of the law of the driving noise process 
for a strong solution of equation \eqref{eq: Tsirel eq} to exist. 
Hirayama--Yano \cite{HY} studied Tsirelson's equation in discrete time 
\begin{align}
X_k = N_k X_{k-1} 
, \quad k \in -\bN , 
\label{eq: Tsirel eq on G}
\end{align}
for processes $ (X_k)_{k \in -\bN} $ and $ (N_k)_{k \in -\bN} $ taking values in a compact group 
and obtained a necessary and sufficient condition 
in terms of infinite product of the driving noise process 
for a strong solution of equation \eqref{eq: Tsirel eq on G} to exist. 
Hirayama--Yano \cite{HYtsirel} studied equation \eqref{eq: Tsirel eq on G} 
for processes taking values in a compact space with semigroup action 
and obtained some sufficient conditions 
for a strong solution of equation \eqref{eq: Tsirel eq on G} to exist 
and not to exist. 

For other contributions of Tsirelson's equation in discrete time, 
see Akahori--Uenishi--Yano \cite{MR2365485} and Takahashi \cite{MR2582432}. 
Several reviews on this topic 
can be found in \cite{YT}, \cite{YY} and \cite{HYtsirel}.

{\bf c). Finite-state Markov chain.} 
As we shall see later in Theorem \ref{thm: main1}, 
for a $ \mu $-random walk $ \{ X,N \} $, 
the process $ X $ is a finite-state Markov chain 
which is stationary, irreducible and aperiodic. 
Yano--Yasutomi \cite{yanoyasu} studied its converse and proved the following: 
Any finite-state Markov chain 
which is stationary, irreducible and aperiodic 
can be realized as a $ \mu $-random walk subject to a synchronizing road coloring.

\subsection{Organization of this paper}

The remainder of this paper is organized as follows. 
In Section \ref{sec: prel}, we introduce some more notations 
and discuss existence and uniqueness of $ \mu $-random walks. 
In Section \ref{sec: examples}, we give several examples 
which help the reader to understand our main theorems deeply. 
Section \ref{sec: main2} is devoted to the proof of Theorem \ref{thm: main intro}. 
In Section \ref{sec: periodic}, we discuss periodic case.

\section{Notations and preliminary facts} \label{sec: prel}

\subsection{Directed graphs and their road colorings}

If the set $ V $ consists of $ m $ elements, 
we may and do write $ V = \{ 1,\ldots,m \} $. 
We shall idenfity $ i \in V $ with $ v_i $ 
where $ \{ v_1,\ldots,v_m \} $ is the standard basis of $ \bR^m $ defined as 
\begin{align}
v_1 = \bmat{1 \\ 0 \\ \vdots \\ 0} , \ 
v_2 = \bmat{0 \\ 1 \\ \vdots \\ 0} , \ldots, \ 
v_m = \bmat{0 \\ 0 \\ \vdots \\ 1} . 
\label{}
\end{align}
We remark that 
the product $ \sigma x $ has two meanings: 
one is the image of the site $ x \in V $ by the mapping $ \sigma \in \Sigma $, 
and the other is 
the usual product among an $ m \times m $-matrix $ [\sigma(y,x)]_{y,x \in V} $ 
and a $ m $-vector (or $ m \times 1 $-matrix) $ x $ in $ \bR^m $. 

The identity mapping $ e $ is identified with the identity matrix. 
We write 
\begin{align}
o_1 = \bmat{1 & 1 & \cdots & 1 \\ 0 & 0 & \cdots & 0 \\ \vdots & \vdots & \ddots & \vdots \\ 0 & 0 & \cdots & 0} 
, \ 
o_2 = \bmat{0 & 0 & \cdots & 0 \\ 1 & 1 & \cdots & 1 \\ \vdots & \vdots & \ddots & \vdots \\ 0 & 0 & \cdots & 0} 
, \ldots, \ 
o_m = \bmat{0 & 0 & \cdots & 0 \\ 0 & 0 & \cdots & 0 \\ \vdots & \vdots & \ddots & \vdots \\ 1 & 1 & \cdots & 1} 
. 
\label{}
\end{align}
The set $ \Sigma $ is a subsemigroup of 
the semigroup consisting of all $ m \times m $ matrices, 
but it is not a group, 
because the elements $ o_1,\ldots,o_m $ 
do not possess its inverse in $ \Sigma $; in fact, we have 
\begin{align}
o_i x = o_i 
, \quad x \in V , \ i=1,\ldots,m . 
\label{}
\end{align}

Let $ \Sigma_0 $ be a subset of $ \Sigma $. 
A sequence $ s = (\sigma_p,\ldots,\sigma_1) $ of $ \Sigma_0 $ 
is called a {\em word in $ \Sigma_0 $}, 
and then $ \abra{s} $ denotes the product $ \sigma_p \cdots \sigma_1 \in \Sigma $. 
We note that the road coloring $ \Sigma_0 $ is synchronizing 
if and only if $ \abra{s} = o_i $ for some word $ s $ in $ \Sigma_0 $ and some $ i=1,\ldots,m $. 
To study the non-synchronizing cases, we introduce the following. 

\begin{Def}
Let $ \Sigma_0 $ be a subset of $ \Sigma $. 
Let $ V_0 $ be a subset of $ V $. 
\subitem (i) 
$ V_0 $ is called {\em synchronizing} 
if $ \abra{s} V_0 $ is a singleton for some word $ s $ in $ \Sigma_0 $. 
\subitem (ii) 
$ V_0 $ is called a {\em deadlock} 
if $ V_0 $ has no synchronizing pair. 
\subitem (iii) 
$ V_0 $ is called {\em stable} 
if the subset $ \abra{s} V_0 $ is synchronizing 
for all word $ s $ in $ \Sigma_0 $. 
\end{Def}

Note that there may exist a pair which is synchronizing but non-stable; 
an example will be given in Section \ref{sec: ex3}. 

\begin{Def}
A subset $ V_0 $ of $ V $ is called an {\em F-clique} if 
$ V_0 $ is a deadlock and is of the form 
$ V_0 = \abra{s} V $ for some word $ s $ in $ \Sigma_0 $. 
\end{Def}

We write $ \sharp(\cdot) $ for the cardinality 
and set 
\begin{align}
\hat{m} = \min \{ \sharp (\abra{s} V) : \text{$ s $ is a word in $ \Sigma_0 $} \} . 
\label{eq: m}
\end{align}
The following lemma can be found in Friedman \cite{MR953004}. 

\begin{Lem}[Friedman \cite{MR953004}] \label{thm: F-clique}
The following assertions hold: 
\subitem {\rm (i)} 
If $ V_0 = \abra{s} V $ for some word $ s $ in $ \Sigma_0 $ 
such that $ \hat{m} = \sharp (\abra{s} V) $, 
then $ V_0 $ is an F-clique. 
\subitem {\rm (ii)} 
For any F-clique $ V_0 $ and for all word $ s $ in $ \Sigma_0 $, 
the subset $ \abra{s} V_0 $ is also an F-clique 
and satisfies $ \sharp(V_0)=\sharp(\abra{s} V_0) $. 
\subitem {\rm (iii)} 
Every F-clique has $ \hat{m} $ elements. 
\end{Lem}

Let us give the proof of Lemma \ref{thm: F-clique} for completeness of this paper. 

\Proof{
(i) Suppose that $ \abra{s} V $ had a synchronizing pair. 
Then there would exist another word $ s' $ in $ \Sigma_0 $ such that 
$ \sharp (\abra{s' s} V) < \sharp (\abra{s} V) = \hat{m} $. 
This contradicts the minimality of $ \hat{m} $. 
Thus we obtain (i). 

(ii) 
Let $ V_0 $ be an F-clique and $ s $ be a word in $ \Sigma_0 $. 
Suppose $ \abra{s} V_0 $ were not an F-clique. 
Then it would admit a synchronizing pair, 
and so would $ V_0 $. This is a contradiction. 

Suppose that $ \sharp(V_0) > \sharp(\abra{s} V_0) $. 
Then $ V_0 $ would admit a synchronizing pair, which is again a contradiction. 

Hence we obtain (ii). 

(iii) 
Let $ V_0 $ be an F-clique 
and $ s $ be a word in $ \Sigma_0 $ such that $ \hat{m} = \sharp(\abra{s} V) $. 
Then, by (ii), we see that  
\begin{align}
\hat{m} \le \sharp(V_0) = \sharp(\abra{s} V_0) \le \sharp(\abra{s} V) = \hat{m}, 
\label{}
\end{align}
which shows that $ \hat{m} = \sharp (V_0) $. 
}

\subsection{Random walks}

We deal with a pair of two processes $ X=(X_k)_{k \in \bZ} $ and $ N=(N_k)_{k \in \bZ} $ 
defined on a common probability space. 
We need the $ \sigma $-fields generated by these processes 
up to time $ k \in \bZ $ given by 
\begin{align}
\cF^X_k = \sigma(X_j : j \le k) 
, \quad 
\cF^N_k = \sigma(N_j : j \le k) 
, \quad 
\cF^{X,N}_k = \sigma(X_j,N_j : j \le k) . 
\label{}
\end{align}
We also need the following $ \sigma $-fields 
for each $ k,l \in \bZ $ with $ k>l $: 
\begin{align}
\cF^X_{k,l} = \sigma(X_j : j=k,k-1,\ldots,l+1) 
, \quad 
\cF^N_{k,l} = \sigma(N_j : j=k,k-1,\ldots,l+1) 
\label{}
\end{align}
and 
\begin{align}
\cF^{X,N}_{k,l} = \sigma(X_j,N_j : j=k,k-1,\ldots,l+1) . 
\label{}
\end{align}

\begin{Def}
Two pairs of processes $ \{ X,N \} $ and $ \{ X',N' \} $ are called {\em identical in law} if 
\begin{align}
\rbra{ (X_{k})_{k \in \bZ},(N_{k})_{k \in \bZ} } 
\dist 
\rbra{ (X'_{k})_{k \in \bZ},(N'_{k})_{k \in \bZ} } . 
\label{}
\end{align}
In this case, we write $ \{ X,N \} \dist \{ X',N' \} $. 
\end{Def}

\begin{Def}
A pair of processes $ \{ X,N \} $ 
is called {\em stationary} if, for any $ n \in \bZ $, 
$ \{ X,N \} $ and $ \{ (X_{k+n})_{k \in \bZ},(N_{k+n})_{k \in \bZ} \} $ 
are identical in law. 
\end{Def}

\begin{Lem} \label{lem: uni}
Let $ \{ X,N \} $ and $ \{ X',N' \} $ be two pairs of processes. 
Then the following assertions hold: 
\subitem {\rm (i)} 
If $ \{ X,N \} \dist \{ X',N' \} $ and if $ \{ X,N \} $ 
is a $ \mu $-random walk (resp. stationary), 
then $ \{ X',N' \} $ is also a $ \mu $-random walk (resp. stationary). 
\subitem {\rm (ii)} 
If $ \{ X,N \} $ and $ \{ X',N' \} $ are $ \mu $-random walks, 
and if $ X_k \dist X'_k $ for all $ k \in \bZ $, then $ \{ X,N \} \dist \{ X',N' \} $. 
\end{Lem}

\Proof{
Claim (i) is obvious. Let us prove Claim (ii). 
Let $ k_0 \in \bZ $. Since $ X_{k_0-1} \dist X'_{k_0-1} $ and $ N \dist N' $ and since 
$ X_k = N_k N_{k-1} \cdots N_{k_0} X_{k_0-1} $ for $ k \ge k_0 $, we see that 
$ (X_k,N_k)_{k \ge k_0} \dist (X'_k,N'_k)_{k \ge k_0} $. 
Since $ k_0 \in \bZ $ is arbitrary, we obtain 
$ (X_k,N_k)_{k \in \bZ} \dist (X'_k,N'_k)_{k \in \bZ} $. 
}

The {\em convolution} of two probability laws $ \mu_1 $ and $ \mu_2 $ on $ \Sigma $ 
will be denoted by $ \mu_1 * \mu_2 $, 
which is a probability law on $ \Sigma $ such that 
\begin{align}
(\mu_1 * \mu_2) (\sigma) 
= \sum_{\sigma_1 \in \Sigma} \mu_1(\sigma_1) 
\sum_{\sigma_2 \in \Sigma} \mu_2(\sigma_2) 
1_{\{ \sigma_1 \sigma_2 = \sigma \}} 
, \quad \sigma \in \Sigma . 
\label{}
\end{align}
For a probability law $ \mu $ on $ \Sigma $ and $ \lambda $ on $ V $, 
the convolution of $ \mu $ and $ \lambda $ 
will also be denoted by $ \mu * \lambda $, 
which is a probability law on $ V $ such that 
\begin{align}
(\mu * \lambda) (y) 
= \sum_{\sigma \in \Sigma} \mu(\sigma) 
\sum_{x \in V} \lambda(x) 
1_{\{ \sigma x = y \}} 
, \quad y \in V . 
\label{}
\end{align}
By \eqref{eq: action}, we see that 
\begin{align}
(\mu_1 * \mu_2) * \lambda = \mu_1 * (\mu_2 * \lambda) . 
\label{}
\end{align}
We write $ \mu^{*1}=\mu $ and define $ \mu^{*n} $ for $ n \ge 2 $ 
recursively by $ \mu^{*n} = \mu^{*(n-1)} * \mu $. 

\begin{Lem} \label{lem: conv eq}
Let $ \{ X,N \} $ be a $ \mu $-random walk. 
For $ k \in \bZ $, let $ \lambda_k $ denote the law of $ X_k $. 
Then the following convolution equation holds: 
\begin{align}
\lambda_k = \mu * \lambda_{k-1} 
, \quad k \in \bZ . 
\label{eq: conv eq}
\end{align}
Conversely, if probability laws $ \mu $ and $ \{ \lambda_k:k \in \bZ \} $ are given 
and the convolution equation \eqref{eq: conv eq} is satisfied, 
then there exists a $ \mu $-random walk $ \{ X,N \} $ 
and $ X_k $ has law $ \lambda_k $ for all $ k \in \bZ $. 
The $ \mu $-random walk is unique up to identity in law. 
\end{Lem}

\Proof{
Let $ \{ X,N \} $ be a $ \mu $-random walk. 
Since $ X_k = N_k X_{k-1} $ and since $ N_k $ is independent of $ X_{k-1} $, 
we obtain \eqref{eq: conv eq}. 

Let $ \mu $ and $ \{ \lambda_k:k \in \bZ \} $ be given 
such that \eqref{eq: conv eq} holds. 
Then, by the Kolmogorov extension theorem, 
we may construct a (possibly time-inhomogeneous) Markov chain $ (N_k,X_{k-1})_{k \in \bZ} $ 
with state space $ V \times \Sigma $ 
so that the marginal law for each $ k \in \bZ $ is given as 
\begin{align}
P(N_k=\sigma, X_{k-1}=x) = \mu(\sigma) \lambda_{k-1}(x) 
, \quad \sigma \in \Sigma , \ x \in V 
\label{}
\end{align}
and the one-step transition probabilities are given as 
\begin{align}
P(N_k=\varsigma, X_{k-1}=y \mid N_{k-1} = \sigma , X_{k-2}=x) 
= \mu(\varsigma) 1_{\{ y=\sigma x \}} 
\label{}
\end{align}
for all $ \sigma,\varsigma \in \Sigma $ and all $ x,y \in V $. 
It is easy to see that the so constructed pair of processes 
$ \{ (X_k)_{k \in \bZ},(N_k)_{k \in \bZ} \} $ is as desired. 
The uniqueness is immediate from Lemma \ref{lem: uni}. 
}

\subsection{Aperiodic case}

Since the index of our $ \mu $-random walk varies in $ \bZ $, 
existence and uniqueness of $ \mu $-random walks are not obvious. 
The following theorem assures the existence and uniqueness 
in the aperiodic case. 

\begin{Thm} \label{thm: main1}
Let $ \mu $ be a probability law on $ \Sigma $ 
and suppose that the directed graph induced by $ \mu $ 
satisfies the assumption {\bf (A)}. 
Then the following assertion holds: 
\subitem {\rm (i)} 
There exists a $ \mu $-random walk $ \{ X,N \} $ in $ V $. 
\subitem {\rm (ii)} 
The $ \mu $-random walk is unique up to identity in law. 
\subitem {\rm (iii)} 
The $ \mu $-random walk $ \{ X,N \} $ is stationary. 
\subitem {\rm (iv)} 
The common law $ \lambda $ of $ X=(X_k)_{k \in \bZ} $ 
is a unique probability law satisfying 
\begin{align}
\mu * \lambda = \lambda . 
\label{eq: conv eq-}
\end{align}
\subitem {\rm (v)} 
The tail $ \sigma $-field $ \cF^{X,N}_{-\infty } := \bigcap_k \cF^{X,N}_k $ 
is a.s. trivial. 
\end{Thm}

This theorem is an immediate consequence of the classical Perron--Frobenius theory 
on infinite product of stochastic matrices. 
We call $ B=(B(y,x))_{x,y=1,\ldots,m} $ a {\em stochastic matrix} if 
$ B(y,x) \ge 0 $ for $ x,y \in V $ and 
\begin{align}
\sum_{y \in V} B(y,x) = 1 
, \quad x \in V . 
\label{}
\end{align}
We call a column vector $ u = {}^{\top} \bmat{u(1) & \cdots & u(m)} $ 
a {\em stochastic vector} if 
$ u(x) \ge 0 $ for all $ x=1,\ldots,m $ and $ \sum_{x=1}^m u(x) = 1 $. 

\begin{Thm}[Perron--Frobenius] \label{thm: PF}
Let $ B=(B(y,x))_{x,y=1,\ldots,m} $ be a stochastic matrix 
and suppose that there exists a positive integer $ r $ 
such that every entry of $ B^r $ is positive. 
Then there exists a stochastic vector $ u $ such that 
\begin{align}
B^n(y,x) \tend{}{n \to \infty } u(y) 
\quad \text{for all $ y = 1,\ldots,m $}, 
\label{}
\end{align}
where $ B^n(y,x) $ is the $ (y,x) $-entry of $ B^n $, i.e., 
$ B^n = (B^n(y,x))_{x,y=1,\ldots,m} $. 
The vector $ u $ is the unique stochastic vector such that 
\begin{align}
B u = u . 
\label{}
\end{align}
\end{Thm}

For the proof of Theorem \ref{thm: PF}, 
see, e.g., \cite{MR2209438}. 

Let us give the proof of Theorem \ref{thm: main1} for completeness of this paper. 

\Proof[Proof of Theorem \ref{thm: main1}]{
(i) 
Define a $ (V \times V) $-matrix $ B=(B(y,x))_{x,y \in V} $ by 
\begin{align}
B = \sum_{\sigma \in \Sigma} \mu(\sigma) \sigma , 
\label{}
\end{align}
or in other words, 
\begin{align}
B(y,x) = \mu( \sigma \in \Sigma : \sigma x = y ) 
, \quad x,y \in V . 
\label{eq: matrix and mu}
\end{align}
Note that $ B $ is a stochastic matrix 
and that $ B(y,x) $ is positive if $ A(y,x) \ge 1 $. 
Since $ (V,A) $ satisfies the assumption {\bf (A)}, 
there exists a positive integer $ r $ 
such that every entry of $ A^r $ is greater than or equal to 1, 
and hence that every entry of $ B^r $ is positive. 
Thus we may apply Theorem \ref{thm: PF} 
to see that 
there exists a probability law $ \lambda $ on $ V $ such that 
\begin{align}
\mu^{*n}(\{ \sigma \in \Sigma : \sigma x = y \}) 
\tend{}{n \to \infty } 
\lambda(y) 
, \quad x,y \in V , 
\label{}
\end{align}
and that $ \lambda $ is the unique probablity law such that 
\begin{align}
\mu * \lambda = \lambda . 
\label{}
\end{align}
By this convolution equation, 
we may apply Lemma \ref{lem: conv eq} to construct $ \{ X,N \} $ such that 
$ N $ has common law $ \mu $ and $ X $ has common law $ \lambda $. 
This is as desired. 

(ii) Let $ \{ X,N \} $ be a $ \mu $-random walk. 
For each $ k \in \bN $, let $ \lambda_k $ denote the law of $ X_k $. 
Then we have 
\begin{align}
\lambda_k = \mu^{*(k-l)} * \lambda_l 
, \quad k,l \in \bZ , \ k>l . 
\label{}
\end{align}
Let $ k \in \bZ $ be fixed. 
Since $ V $ is finite, there exist a subsequence $ l(n) \to -\infty $ 
and a probability law $ \tilde{\lambda} $ such that 
$ \lambda_{l(n)} \tend{\rm w}{} \tilde{\lambda} $. 
Then, for any $ y \in V $, we have 
\begin{align}
\lambda_k(y) 
=& \mu^{*(k-l(n))} * \lambda_{l(n)}(y) 
\label{} \\
=& \sum_{x \in V} \mu^{*(k-l(n))}(\{ \sigma \in \Sigma : \sigma x = y \}) \lambda_{l(n)}(x) 
\label{} \\
\tend{}{n \to \infty }& 
\lambda(y) \sum_{x \in V} \tilde{\lambda}(x) 
= \lambda(y) . 
\label{}
\end{align}
This shows that $ X_k $ has law $ \lambda $ for each $ k \in \bZ $. 
This proves uniqueness by Lemma \ref{lem: uni}. 

Claims (iii) and (iv) have already been proved. 

(v) Suppose that 
there were $ A \in \cF^{X,N}_{-\infty } $ such that $ 0<P(A)<1 $. 
Define $ P'= P(\cdot \mid A) $. 
Then it is easy to see that 
$ \{ X,N \} $ under $ P' $ is also a $ \mu $-random walk. 
Thus Claim (ii) shows that $ P'((X,N) \in \cdot) = P((X,N) \in \cdot) $. 
This contradicts the fact that $ P'(A) = 1 > P(A) $. 
This proves Claim (v). 

The proof is now complete. 
}

\subsection{Proofs of easy parts of Theorem \ref{thm: main intro}}

Let us prove easy parts of Theorem \ref{thm: main intro}. 

\Proof[Proof of {\rm [(i) $ \Rightarrow $ (ii)]} of Theorem \ref{thm: main intro}.]{
Suppose that (i) holds. 
Let $ s = (\sigma_p,\ldots,\sigma_1) $ be a word in $ \Supp(\mu) $ 
such that $ \abra{s} V $ is a singleton. 
Define, for each $ k \in \bZ $, a random time 
\begin{align}
T(k) = \max \{ l = k-1,k-2,\ldots : \ 
N_{l+p} = \sigma_p, \ldots, N_{l+1} = \sigma_1 \} . 
\label{}
\end{align}
Here we note that the random time $ k-T(k) $ is a stopping time 
with respect to the filtration $ \{ \cF^N_{k,k-n} : n=1,2,\ldots \} $. 
Then we see that $ T(k) $ is finite a.s. and that 
\begin{align}
\lim_{l \to -\infty } N_{k,l} = N_{k,T(k)+p} 
\quad \text{a.s.} 
\label{}
\end{align}
Thus we obtain (ii). 
}

\Proof[Proof of {\rm [(ii) $ \Rightarrow $ (iii)]} of Theorem \ref{thm: main intro}.]{
Suppose that (ii) holds. Denote 
\begin{align}
Y_k = \lim_{l \to -\infty } N_k N_{k-1} \cdots N_{l+1} 
\quad \text{a.s. for $ k \in \bZ $}. 
\label{}
\end{align}
Then, for any fixed $ x_0 \in V $, 
a pair of processes 
\begin{align}
\{ (Y_k x_0)_{k \in \bZ}, (N_k)_{k \in \bZ} \} 
\label{}
\end{align}
is a strong $ \mu $-random walk which is identical in law to $ \{ X,N \} $. 
This proves (iii). 
}

\section{Illustrative examples} \label{sec: examples}

Before proceeding to prove our main theorems, 
we give illustrative examples. 

\

\minip{7cm}{
\unitlength 0.1in
\begin{picture}( 20.0000, 16.0000)(  0.0000,-16.0000)
%
\special{pn 8}%
\special{ar 1000 200 200 200  0.0000000 6.2831853}%
%
\special{pn 8}%
\special{ar 200 1400 200 200  0.0000000 6.2831853}%
%
\special{pn 8}%
\special{ar 1800 1400 200 200  0.0000000 6.2831853}%
%
\special{pn 8}%
\special{pa 810 270}%
\special{pa 210 1190}%
\special{fp}%
\special{pa 960 400}%
\special{pa 380 1280}%
\special{fp}%
%
\special{pn 8}%
\special{pa 380 1330}%
\special{pa 1600 1330}%
\special{fp}%
\special{pa 380 1490}%
\special{pa 1610 1490}%
\special{fp}%
%
\special{pn 8}%
\special{pa 1070 390}%
\special{pa 1650 1250}%
\special{fp}%
\special{pa 1190 260}%
\special{pa 1810 1200}%
\special{fp}%
%
\special{pn 8}%
\special{pa 200 1060}%
\special{pa 200 1200}%
\special{fp}%
\special{pa 200 1200}%
\special{pa 350 1130}%
\special{fp}%
\special{pa 820 460}%
\special{pa 950 400}%
\special{fp}%
\special{pa 950 390}%
\special{pa 970 560}%
\special{fp}%
\special{pa 510 1240}%
\special{pa 390 1320}%
\special{fp}%
\special{pa 390 1320}%
\special{pa 520 1420}%
\special{fp}%
%
\special{pn 20}%
\special{pa 380 1490}%
\special{pa 1610 1490}%
\special{fp}%
\special{pa 1480 1410}%
\special{pa 1610 1490}%
\special{fp}%
\special{pa 1610 1490}%
\special{pa 1460 1590}%
\special{fp}%
\special{pa 1650 1250}%
\special{pa 1070 390}%
\special{fp}%
\special{pa 1060 580}%
\special{pa 1080 420}%
\special{fp}%
\special{pa 1080 430}%
\special{pa 1230 490}%
\special{fp}%
\special{pa 1190 270}%
\special{pa 1810 1190}%
\special{fp}%
\special{pa 1670 1140}%
\special{pa 1820 1210}%
\special{fp}%
\special{pa 1810 1200}%
\special{pa 1840 1050}%
\special{fp}%
\put(10.0000,-2.2000){\makebox(0,0){1}}%
\put(1.9000,-14.1000){\makebox(0,0){2}}%
\put(18.0000,-14.0000){\makebox(0,0){3}}%
\end{picture}%
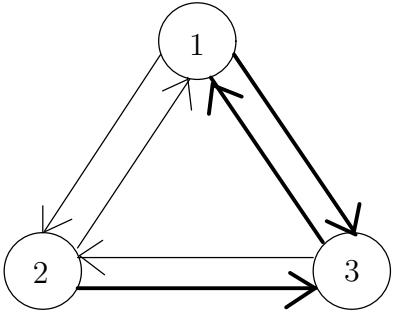
\\ Figure 3. 
}
\minip{7cm}{
\unitlength 0.1in
\begin{picture}( 20.0000, 16.0000)(  0.0000,-16.0000)
%
\special{pn 8}%
\special{ar 1000 200 200 200  0.0000000 6.2831853}%
%
\special{pn 8}%
\special{ar 200 1400 200 200  0.0000000 6.2831853}%
%
\special{pn 8}%
\special{ar 1800 1400 200 200  0.0000000 6.2831853}%
%
\special{pn 8}%
\special{pa 810 270}%
\special{pa 210 1190}%
\special{fp}%
\special{pa 960 400}%
\special{pa 380 1280}%
\special{fp}%
%
\special{pn 8}%
\special{pa 380 1330}%
\special{pa 1600 1330}%
\special{fp}%
\special{pa 380 1490}%
\special{pa 1610 1490}%
\special{fp}%
%
\special{pn 8}%
\special{pa 1070 390}%
\special{pa 1650 1250}%
\special{fp}%
\special{pa 1190 260}%
\special{pa 1810 1200}%
\special{fp}%
%
\special{pn 8}%
\special{pa 240 980}%
\special{pa 210 1200}%
\special{fp}%
\special{pa 200 1190}%
\special{pa 360 1090}%
\special{fp}%
\special{pa 1440 1380}%
\special{pa 1600 1490}%
\special{fp}%
\special{pa 1590 1490}%
\special{pa 1450 1600}%
\special{fp}%
\special{pa 1190 450}%
\special{pa 1200 260}%
\special{fp}%
\special{pa 1190 260}%
\special{pa 1360 340}%
\special{fp}%
%
\special{pn 20}%
\special{pa 960 400}%
\special{pa 370 1290}%
\special{fp}%
\special{pa 390 1330}%
\special{pa 1600 1330}%
\special{fp}%
\special{pa 1650 1250}%
\special{pa 1080 400}%
\special{fp}%
\special{pa 780 480}%
\special{pa 950 420}%
\special{fp}%
\special{pa 970 410}%
\special{pa 930 620}%
\special{fp}%
\special{pa 1470 1230}%
\special{pa 1650 1250}%
\special{fp}%
\special{pa 1650 1250}%
\special{pa 1670 1100}%
\special{fp}%
\special{pa 520 1210}%
\special{pa 390 1310}%
\special{fp}%
\special{pa 390 1310}%
\special{pa 540 1430}%
\special{fp}%
\put(10.0000,-2.1000){\makebox(0,0){1}}%
\put(2.0000,-14.1000){\makebox(0,0){2}}%
\put(18.0000,-14.0000){\makebox(0,0){3}}%
\end{picture}%
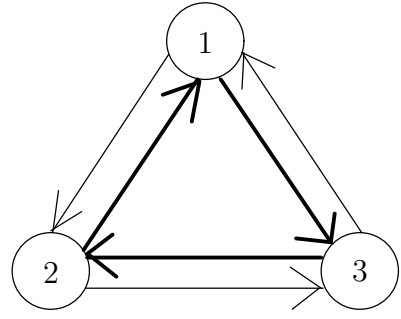
\\ Figure 4. 
}

\subsection{Synchronizing case} \label{sec: ex1}

Let $ V=\{ 1,2,3 \} $ and consider the following adjacency matrix: 
\begin{align}
A = 
\bmat{
A(1,1) & A(1,2) & A(1,3) \\
A(2,1) & A(2,2) & A(2,3) \\
A(3,1) & A(3,2) & A(3,3) \\
} 
= 
\bmat{
0 & 1 & 1 \\
1 & 0 & 1 \\
1 & 1 & 0 
} . 
\label{eq: ex1}
\end{align}
We can easily verify that the graph $ (V,A) $ satisfies the assumption {\bf (A)}. 
Consider a road coloring $ \{ \sigma^{(1)},\sigma^{(2)} \} $ of $ (V,A) $ given as 
\begin{align}
\sigma^{(1)} = 
\bmat{
0 & 0 & 1 \\
0 & 0 & 0 \\
1 & 1 & 0 
} , \quad 
\sigma^{(2)} = 
\bmat{
0 & 1 & 0 \\
1 & 0 & 1 \\
0 & 0 & 0 
} . 
\label{}
\end{align}
See Figure 3, where the thick roads are colored $ \sigma^{(1)} $ 
and the thin ones $ \sigma^{(2)} $. 
Since 
\begin{align}
\sigma^{(1)} \sigma^{(2)} = o_3 , 
\label{}
\end{align}
we see that the road coloring $ \{ \sigma^{(1)},\sigma^{(2)} \} $ 
is synchronizing. 

Let $ p,q > 0 $ with $ p+q=1 $ 
and let $ \mu $ be a probability law on $ \Sigma $ such that 
\begin{align}
\mu(\{ \sigma^{(1)} \}) = p 
, \quad 
\mu(\{ \sigma^{(2)} \}) = q . 
\label{eq: def of mu}
\end{align}
Let $ \{ X,N \} $ be a $ \mu $-random walk in $ (V,A) $. 
Solving equation \eqref{eq: conv eq-}, 
we see that the common law $ \lambda $ of $ X $ is given as 
\begin{align}
\bmat{\lambda(1) \\ \lambda(2) \\ \lambda(3) } 
= \frac{1}{2+pq} \bmat{1-pq \\ q+pq \\ p+pq} 
\label{}
\end{align}
where we write $ \lambda(i) $ simply for $ \lambda(\cbra{i}) $, $ i=1,2,3 $. 
For each $ k \in \bZ $, we have 
\begin{align}
T(k) = \inf \{ n=1,2,\ldots, : \ N_{k-n+1} = \sigma^{(1)} , \ N_{k-n+2} = \sigma^{(2)} \} 
\label{}
\end{align}
and then we have 
\begin{align}
X_k = N_{k,k-T(k)+2} v_3 
, \quad k \in \bZ , 
\label{}
\end{align}
which shows that the $ \mu $-random walk $ \{ X,N \} $ is strong.

\subsection{Non-synchronizing case: an easy example} \label{sec: ex2}

Let $ V=\{ 1,2,3 \} $ and let $ A $ as defined in \eqref{eq: ex1}. 
Consider a road coloring $ \{ \sigma^{(1)},\sigma^{(2)} \} $ of $ (V,A) $ given as 
\begin{align}
\sigma^{(1)} = 
\bmat{
0 & 1 & 0 \\
0 & 0 & 1 \\
1 & 0 & 0 
} , \quad 
\sigma^{(2)} = 
\bmat{
0 & 0 & 1 \\
1 & 0 & 0 \\
0 & 1 & 0 
} . 
\label{eq: ex2-1}
\end{align}
See Figure 4, where the thick roads are colored $ \sigma^{(1)} $ 
and the thin ones $ \sigma^{(2)} $. 
Since 
\begin{align}
\sigma^{(1)} \sigma^{(2)} = 
\sigma^{(2)} \sigma^{(1)} = 
{\rm id} , 
\label{}
\end{align}
the set $ G = \{ \sigma^{(1)},\sigma^{(2)},{\rm id} \} $ is a group, 
and thus we see that 
the road coloring $ \{ \sigma^{(1)},\sigma^{(2)} \} $ is non-synchronizing. 

Let $ \mu $ be a probability law on $ \Sigma $ such that 
$ \Supp(\mu) = \{ \sigma^{(1)},\sigma^{(2)} \} $. 
Let $ \{ X,N \} $ be a $ \mu $-random walk in $ (V,A) $. 
Solving equation \eqref{eq: conv eq-}, 
we see that the common law $ \lambda $ of $ X $ is uniform on $ V $, i.e., 
\begin{align}
\bmat{\lambda(1) \\ \lambda(2) \\ \lambda(3) } 
= \frac{1}{3} \bmat{1 \\ 1 \\ 1} . 
\label{}
\end{align}
Now, by Proposition \ref{prop: typical}, 
we conclude that $ \{ X,N \} $ is non-strong.

\subsection{Non-synchronizing case: a difficult example} \label{sec: ex3}

\begin{center}
\unitlength 0.1in
\begin{picture}( 20.0000, 18.0000)(  0.0000,-18.0000)
%
\special{pn 8}%
\special{ar 200 800 200 200  0.0000000 6.2831853}%
%
\special{pn 8}%
\special{ar 1000 200 200 200  0.0000000 6.2831853}%
%
\special{pn 8}%
\special{ar 1800 800 200 200  0.0000000 6.2831853}%
%
\special{pn 8}%
\special{ar 600 1600 200 200  0.0000000 6.2831853}%
%
\special{pn 8}%
\special{ar 1400 1600 200 200  0.0000000 6.2831853}%
\put(10.0000,-2.0000){\makebox(0,0){1}}%
\put(2.0000,-8.0000){\makebox(0,0){2}}%
\put(6.0000,-16.0000){\makebox(0,0){3}}%
\put(14.0000,-16.0000){\makebox(0,0){4}}%
\put(18.0000,-8.0000){\makebox(0,0){5}}%
%
\special{pn 20}%
\special{ar 1800 550 184 184  2.4858970 6.2831853}%
\special{ar 1800 550 184 184  0.0000000 0.6202495}%
%
\special{pn 20}%
\special{pa 1550 610}%
\special{pa 1650 660}%
\special{fp}%
\special{pa 1640 640}%
\special{pa 1700 550}%
\special{fp}%
\special{pa 880 370}%
\special{pa 400 740}%
\special{fp}%
\special{pa 410 620}%
\special{pa 400 740}%
\special{fp}%
\special{pa 400 740}%
\special{pa 530 740}%
\special{fp}%
\special{pa 290 980}%
\special{pa 500 1410}%
\special{fp}%
\special{pa 380 1380}%
\special{pa 490 1410}%
\special{fp}%
\special{pa 490 1410}%
\special{pa 530 1290}%
\special{fp}%
\special{pa 800 1600}%
\special{pa 1180 1600}%
\special{fp}%
\special{pa 1100 1530}%
\special{pa 1200 1600}%
\special{fp}%
\special{pa 1200 1600}%
\special{pa 1120 1670}%
\special{fp}%
\special{pa 1380 1410}%
\special{pa 1090 390}%
\special{fp}%
\special{pa 1060 540}%
\special{pa 1090 410}%
\special{fp}%
\special{pa 1090 410}%
\special{pa 1220 500}%
\special{fp}%
%
\special{pn 8}%
\special{pa 800 260}%
\special{pa 310 610}%
\special{fp}%
\special{pa 310 510}%
\special{pa 300 610}%
\special{fp}%
\special{pa 300 610}%
\special{pa 410 610}%
\special{fp}%
\special{pa 410 810}%
\special{pa 1600 810}%
\special{fp}%
\special{pa 1490 740}%
\special{pa 1600 810}%
\special{fp}%
\special{pa 1600 810}%
\special{pa 1480 900}%
\special{fp}%
\special{pa 1240 1460}%
\special{pa 370 900}%
\special{fp}%
\special{pa 410 1010}%
\special{pa 380 900}%
\special{fp}%
\special{pa 380 900}%
\special{pa 490 900}%
\special{fp}%
\special{pa 750 1460}%
\special{pa 1640 940}%
\special{fp}%
\special{pa 1510 930}%
\special{pa 1630 930}%
\special{fp}%
\special{pa 1630 930}%
\special{pa 1590 1070}%
\special{fp}%
\special{pa 1760 1010}%
\special{pa 1530 1430}%
\special{fp}%
\special{pa 1510 1320}%
\special{pa 1520 1440}%
\special{fp}%
\special{pa 1520 1440}%
\special{pa 1640 1400}%
\special{fp}%
\end{picture}%
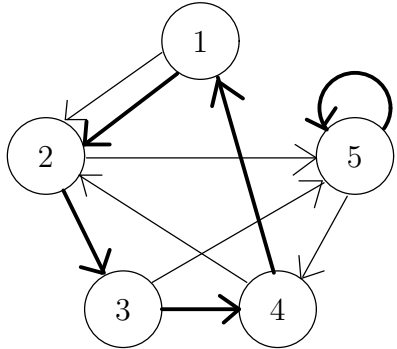
\\ Figure 5. 
\end{center}

Let $ V=\{ 1,2,3,4,5 \} $ 
and consider the following adjacency matrix: 
\begin{align}
A = 
\bmat{
0 & 0 & 0 & 1 & 0 \\
2 & 0 & 0 & 1 & 0 \\
0 & 1 & 0 & 0 & 0 \\
0 & 0 & 1 & 0 & 1 \\
0 & 1 & 1 & 0 & 1 
} . 
\label{eq: ex3}
\end{align}
We can easily verify that the graph $ (V,A) $ satisfies the assumption {\bf (A)}. 
Consider a road coloring $ \{ \sigma^{(1)},\sigma^{(2)} \} $ of $ (V,A) $ given as 
\begin{align}
\sigma^{(1)} = 
\bmat{
0 & 0 & 0 & 1 & 0 \\
1 & 0 & 0 & 0 & 0 \\
0 & 1 & 0 & 0 & 0 \\
0 & 0 & 1 & 0 & 0 \\
0 & 0 & 0 & 0 & 1 
} , \quad 
\sigma^{(2)} = 
\bmat{
0 & 0 & 0 & 0 & 0 \\
1 & 0 & 0 & 1 & 0 \\
0 & 0 & 0 & 0 & 0 \\
0 & 0 & 0 & 0 & 1 \\
0 & 1 & 1 & 0 & 0 
} , \quad 
\label{}
\end{align}
See Figure 5, where 
the thick roads are colored $ \sigma^{(1)} $ 
and the thin ones $ \sigma^{(2)} $. 
It is easy to see that the F-cliques are 
\begin{align}
\{ 1,3,5 \} \quad \text{and} \quad \{ 2,4,5 \} . 
\label{}
\end{align}
In particular, the road coloring is non-synchronizing. 
The pairs 
\begin{align}
\{ 1,2 \}, \{ 1,4 \}, \{ 2,3 \}, \{ 3,4 \} 
\label{eq: non-stable pairs}
\end{align}
are all synchronizing. Since 
\begin{align}
\{ 1,4 \} 
\tend{\sigma^{(1)}}{} 
\{ 1,2 \} 
\tend{\sigma^{(2)}}{} 
\{ 2,5 \} 
, \quad 
\{ 2,3 \} 
\tend{\sigma^{(1)}}{} 
\{ 3,4 \} 
\tend{\sigma^{(2)}}{} 
\{ 2,5 \} , 
\label{}
\end{align}
and since $ \{ 2,5 \} $ is a deadlock, 
we see that the pairs \eqref{eq: non-stable pairs} are non-stable. 

Let $ p,q > 0 $ with $ p+q=1 $ 
and let $ \mu $ be as defined in \eqref{eq: def of mu}. 
Let $ \{ X,N \} $ be a $ \mu $-random walk in $ (V,A) $. 
Solving equation \eqref{eq: conv eq-}, 
we see that the common law $ \lambda $ of $ X $ is given as 
\begin{align}
\bmat{\lambda(1) \\ \lambda(2) \\ \lambda(3) \\ \lambda(4) \\ \lambda(5) } 
= \frac{1}{3(1+p)} \bmat{p \\ 1 \\ p \\ 1 \\ 1+p} . 
\label{}
\end{align}
Thus we may find the uniformity: 
\begin{align}
\lambda(\{ 1,2 \}) = \lambda(\{ 3,4 \}) = \lambda(\{ 5 \}) . 
\label{}
\end{align}
Note that 
this is a special case of Theorem \ref{thm: main3} given in the next section. 
One may expect that some symmetry lies behind this uniformity, 
but it seems hidden because $ \{ 1,2 \} $ and $ \{ 3,4 \} $ cannot be interchanged. 
We will reveal a certain hidden symmetry behind this uniformity 
in the proof of Claim \eqref{eq: main2}.

\section{Non-strongness of the $ \mu $-random walk in the non-synchronizing case} \label{sec: main2}

This section is devoted to the proof of Claim \eqref{eq: main2}, 
which will complete the proof of Theorem \ref{thm: main intro}. 

\subsection{Uniformity}

Suppose that $ \{ X,N \} $ and $ \lambda $ be as in Theorem \ref{thm: main1}. 
Let $ \Supp(\mu) = \{ \sigma^{(1)},\ldots,\sigma^{(d)} \} $ 
and suppose that $ \Supp(\mu) $ is non-synchronizing. 
Then, by (i) of Lemma \ref{thm: F-clique}, 
there exists a word $ s = (\sigma_p,\ldots,\sigma_1) $ in $ \Supp(\mu) $ 
such that $ \abra{s} V $ is an F-clique. 
We enumerate $ \abra{s} V $ as 
\begin{align}
\abra{s} V = \{ \hat{x}_1,\ldots,\hat{x}_{\hat{m}} \} 
\label{}
\end{align}
where 
\begin{align}
\hat{m} = \min \{ \sharp (\abra{s} V) : \text{$ s $ is a word in $ \Supp(\mu) $} \} . 
\label{}
\end{align}
Since $ \Supp(\mu) $ is non-synchronizing, we have 
\begin{align}
\hat{m} \ge 2 . 
\label{}
\end{align}
Set 
\begin{align}
V_i = \{ x \in V : \ \abra{s} x = \hat{x}_i \} 
, \quad i = 1, \ldots, \hat{m} . 
\label{}
\end{align}
Then the family $ \{ V_1,\ldots,V_{\hat{m}} \} $ is a partition of the state space $ V $. 
Note that this partition of $ V $ may depend on the choice of the word $ s $ in $ \Supp(\mu) $ 
such that $ \abra{s} V $ is an F-clique. 
The following theorem is crucial to our proof of Claim \eqref{eq: main2}, 
which does not matter whatever we choose as such a word $ s $ in $ \Supp(\mu) $. 

\begin{Thm} \label{thm: main3}
It holds that 
\begin{align}
\lambda(V_1) = \cdots = \lambda(V_{\hat{m}}) = \frac{1}{\hat{m}} . 
\label{eq: uniformity}
\end{align}
\end{Thm}

We shall postpone the proof of Theorem \ref{thm: main3} 
until Section \ref{sec: new RW}.

\subsection{Constructing a permutation process}

Denote $ \hat{V} = \{ 1,\ldots,\hat{m} \} $ 
and its permutation group by $ \fS(\hat{V}) $. 
We decompose $ V $ into the disjoint union $ \bigcup_{i=1}^{\hat{m}} V_i $ where 
\begin{align}
V_i = \{ x \in V : \abra{s} x = \hat{x}_i \} 
, \quad i \in \hat{V} . 
\label{}
\end{align}
By (ii) and (iii) of Lemma \ref{thm: F-clique}, we see that 
\begin{align}
\sharp (\abra{s} \sigma \abra{s} V) = \hat{m} 
\quad \text{for any} \ \sigma \in \Sigma , 
\label{}
\end{align}
and hence 
\begin{align}
(\abra{s} \sigma) \{ \hat{x}_1,\ldots,\hat{x}_{\hat{m}} \} 
= \{ \hat{x}_1,\ldots,\hat{x}_{\hat{m}} \} . 
\label{}
\end{align}
This yields that there exists a mapping 
$ M:\Sigma \ni \sigma \mapsto M[\sigma] \in \fS(\hat{V}) $ 
such that 
\begin{align}
M[\sigma](i) = j 
\quad \text{if and only if} \quad 
\sigma \hat{x}_i \in V_j , 
\label{}
\end{align}
where $ i,j \in \hat{V} $. 

Let $ \{ X,N \} $ be a $ \mu $-random walk in $ (V,A) $. 
Set $ T_1=0 $ and define $ T_0,T_{-1},\ldots $ recursively by 
\begin{align}
T_{\kappa} = \max \{ l \le T_{\kappa+1}-p : N_{l+p} = \sigma_p , 
\ldots, N_{l+2} = \sigma_2, N_{l+1} = \sigma_1\} . 
\label{}
\end{align}
By the second Borel--Cantelli lemma, 
we see that the decreasing sequence $ (T_{\kappa })_{\kappa \in -\bN} $ is well-defined a.s. 
Note that, for any $ \kappa \in -\bN $ and any $ l \in -\bN $, we have 
\begin{align}
\{ T_{\kappa} = l \} \in \cF^N_{0,l} . 
\label{}
\end{align}
Now we define an $ \fS(\hat{V}) $-valued process $ (\hat{N}_{\kappa })_{\kappa \in -\bN} $ as 
\begin{align}
\hat{N}_{\kappa } = M \sbra{ 
N_{T_{\kappa }} \cdots N_{T_{\kappa -1}+p+2} N_{T_{\kappa -1}+p+1} } 
\quad \text{if $ T_{\kappa -1} < T_{\kappa} - p $}
\label{}
\end{align}
and $ \hat{N}_{\kappa } = $ identity if $ T_{\kappa -1} = T_{\kappa} - p $. 
We may write $ \hat{\mu} $ for the law of $ \hat{N}_0 $ on $ \fS(\hat{V}) $. 
Then it is obvious that $ (\hat{N}_{\kappa })_{\kappa \in -\bN} $ 
has common law $ \hat{\mu} $ 
since 
\begin{align}
N_{T_{\kappa }} \cdots N_{T_{\kappa -1}+p+2} N_{T_{\kappa -1}+p+1} 
\dist N_0 \cdots N_{T_{-1}+p+2} N_{T_{-1}+p+1} . 
\label{}
\end{align}

By the aperiodicity assumption, there exists a constant $ r $ such that 
from any $ x \in V $ to any $ y \in V $ there exists a path of length $ r $. 
For a technical reason, we introduce the following assumption: 
\begin{align}
\begin{split}
p > r \quad & \text{and, for any $ q=p-1,p-2,\ldots,p-r $,} 
\\
& (\sigma_p,\sigma_{p-1},\ldots,\sigma_{p-q+2},\sigma_{p-q+1}) \neq 
(\sigma_q,\sigma_{q-1},\ldots,\sigma_2,\sigma_1) . 
\end{split}
\label{eq: ass}
\end{align}
To prove Theorem \ref{thm: main3}, we may assume \eqref{eq: ass} 
without loss of generality. 
For this, it suffices to replace $ s $ by another word $ \tilde{s} $ in $ \Supp(\mu) $ 
defined as follows. 
Set $ \tilde{p} = p+2r $ and 
\begin{align}
\tilde{\sigma}_i = 
\begin{cases}
\sigma_i & \text{if $ i=1,2,\ldots,p $}, \\
\sigma^{(1)} & \text{if $ i=p+1,p+2,\ldots,p+r $}, \\
\sigma^{(2)} & \text{if $ i=p+r+1,p+r+2,\ldots,p+2r $}, 
\end{cases}
\label{}
\end{align}
and then define 
\begin{align}
\tilde{s} = (\tilde{\sigma}_{\tilde{p}}, \ldots \tilde{\sigma}_2, \tilde{\sigma}_1) . 
\label{}
\end{align}
Then it is obvious that 
the sequence $ \tilde{s} $ satisfies \eqref{eq: ass}. 
By (ii) and (iii) of Lemma \ref{thm: F-clique}, we see that 
$ \abra{\tilde{s}} V $ is also an F-clique, and that 
$ \abra{\tilde{s}} V = \{ \tilde{x}_1,\ldots,\tilde{x}_{\hat{m}} \} $. 
Set 
\begin{align}
\tilde{V}_i = \{ x \in V : \abra{\tilde{s}} x = \tilde{x}_i \} 
, \quad i=1,\ldots,\hat{m} . 
\label{}
\end{align}
We then note that $ (\tilde{V}_1,\ldots,\tilde{V}_{\hat{m}}) $ 
is a permutation of $ (V_1,\ldots,V_{\hat{m}}) $, 
which shows that the replacement of $ s $ by $ \tilde{s} $ 
does not matter in the proof of Theorem \ref{thm: main3}. 

\begin{Lem} \label{lem: lem}
Suppose that \eqref{eq: ass} holds. 
Then, for any $ i,j \in \hat{V} $, it holds that 
\begin{align}
\hat{\mu}^{*n} \rbra{ \cbra{ \pi \in \fS(\hat{V}) : \pi(i)=j } } 
\tend{}{n \to \infty } 
\frac{1}{\hat{m}} . 
\label{}
\end{align}
\end{Lem}

\Proof{
By the strong-connectedness property of $ (V,A) $ 
we see that, for any $ i,j \in \hat{V} $, 
there exists a word $ s' = (\sigma'_r,\ldots,\sigma'_1) $ in $ \Supp(\mu) $ 
such that 
\begin{align}
\abra{s'} \hat{x}_i \in \hat{V}_j . 
\label{}
\end{align}
We define an event $ B $ by 
\begin{align}
\begin{split}
B = \{ & N_0 = \sigma'_r, \ldots, \ 
N_{-r+2} = \sigma'_2 , \ N_{-r+1} = \sigma'_1 , 
\\
& N_{-r} = \sigma_p ,\ldots, \ N_{-r-p+2} = \sigma_2 , \ N_{-r-p+1} = \sigma_1 \} . 
\end{split}
\label{}
\end{align}
Since we assume that \eqref{eq: ass} holds, we see that 
\begin{align}
T_{-1} = -r-p \ \text{and} \ \hat{N}_0 (i) = j 
\ \text{on the event $ B $}. 
\label{}
\end{align}
Therefore, for any $ i,j \in \hat{V} $, we have 
\begin{align}
\hat{\mu}(\{ \pi \in \fS(\hat{V}) : \pi(i)=j \}) \ge P(B) > 0 . 
\label{}
\label{}
\end{align}
Now we may apply Theorem \ref{thm: PF} to see that 
there exists a probability law $ (\rho(j):j \in \hat{V}) $ such that 
\begin{align}
\hat{\mu}^{*n}(\{ \pi \in \fS(\hat{V}) : \pi(i)=j \}) 
\tend{}{n \to \infty } 
\rho(j) 
, \quad i,j \in \hat{V} . 
\label{}
\end{align}
Since, for any $ n $ and $ j \in \hat{V} $, we have 
\begin{align}
\sum_{i \in \hat{V}} \hat{\mu}^{*n}(\{ \pi \in \fS(\hat{V}) : \pi(i)=j \}) 
= \sum_{i \in \hat{V}} \hat{\mu}^{*n}(\{ \pi \in \fS(\hat{V}) : \pi(j)=i \}) = 1 , 
\label{}
\end{align}
we see that $ \rho(j) = 1/\hat{m} $ for all $ j \in \hat{V} $. 
The proof is now complete. 
}

\subsection{Constructing a new random walk} \label{sec: new RW}

Define $ \hat{X}_{\kappa } \in \hat{V} $ as follows: 
\begin{align}
\hat{X}_{\kappa } = i 
\quad \text{if} \quad 
X_{T_{\kappa }} \in \hat{V}_i . 
\label{}
\end{align}
We will prove gradually that the process 
$ \{ \hat{X},\hat{N} \} 
= \{ (\hat{X}_{\kappa })_{\kappa \in - \bN},(\hat{N}_{\kappa })_{\kappa \in - \bN} \} $ 
is a $ \hat{\mu} $-random walk indexed by $ -\bN $. 
The first step is the following.

\begin{Lem} \label{lem: lem2}
Suppose that \eqref{eq: ass} holds. 
Then, for fixed $ \kappa \in -\bN $, the following assertions hold: 
\subitem (i) 
$ \hat{X}_{\kappa } = \hat{N}_{\kappa } \hat{X}_{\kappa -1} $ holds a.s.; 
\subitem (ii) 
$ \hat{X}_{\kappa -1} $ is indepenent of $ \hat{N}_{\kappa } $. 
\subitem (iii) 
$ P(\hat{X}_{\kappa -1}=i) = \lambda(\hat{V}_i) $ for all $ i \in \hat{V} $. 
\end{Lem}

\Proof{
Claim (i) is obvious by definition. 
Let us prove (ii) and (iii) at the same time. 
Let $ l \in -\bN $, $ \pi \in \fS(\hat{V}) $ and $ i \in \hat{V} $. 
Then we have 
\begin{align}
& P \rbra{ \hat{N}_{\kappa } = \pi , \ \hat{X}_{\kappa -1} = i , \ T_{\kappa -1} = l } 
\label{} \\
=& P \rbra{ \hat{N}_{\kappa } = \pi , \ X_l \in V_i , \ T_{\kappa -1} = l } 
\label{} \\
=& P \rbra{ \hat{N}_{\kappa } = \pi , \ T_{\kappa -1} = l} P(X_l=i) 
& \text{(by independence)} 
\label{} \\
=& P \rbra{ \hat{N}_{\kappa } = \pi , \ T_{\kappa -1} = l} \lambda(V_i) 
& \text{(by stationarity).} 
\label{}
\end{align}
Summing up by $ l \in -\bN $, we obtain 
\begin{align}
P \rbra{ \hat{N}_{\kappa } = \pi , \ \hat{X}_{\kappa -1} = i } 
= P \rbra{ \hat{N}_{\kappa } = \pi } \lambda(V_i) , 
\label{eq: lambda Vi}
\end{align}
which proves Claims (ii) and (iii). 
}

The second step is to prove Theorem \ref{thm: main3}. 

\Proof[Proof of Theorem \ref{thm: main3}]{
Define a probability law $ \hat{\lambda} $ on $ \hat{V} $ by 
\begin{align}
\hat{\lambda}(i) = \lambda(\hat{V}_i) 
\quad \text{for all $ i \in \hat{V} $}. 
\label{}
\end{align}
By Lemma \ref{lem: lem2}, we have 
\begin{align}
\hat{\lambda} = \hat{\mu} * \hat{\lambda} . 
\label{}
\end{align}
Iterating this convolution equation, we obtain, for any fixed $ i \in \hat{V} $, 
\begin{align}
\hat{\lambda}(\{ i \}) 
=& (\hat{\mu}^{*n} * \hat{\lambda})(\{ i \}) 
\label{} \\
=& \sum_{j \in \hat{V}} 
\hat{\mu}^{*n}(\{ \pi \in \fS(\hat{V}) : \pi(j)=i \}) \hat{\lambda}(\{ j \}) 
\label{} \\
\tend{}{n \to \infty }& 
\frac{1}{\hat{m}} \sum_{j \in \hat{V}} \hat{\lambda}(\{ j \}) 
= \frac{1}{\hat{m}} 
\label{}
\end{align}
where we have used Lemma \ref{lem: lem}. 
This completes the proof. 
}

The third step is the following, 
which reveals a symmetry hidden behind the uniformity \eqref{eq: uniformity} 
in Theorem \ref{thm: main3}. 

\begin{Thm}
Suppose that \eqref{eq: ass} holds. 
Then $ \{ (\hat{X}_{\kappa })_{\kappa \in -\bN},(\hat{N}_{\kappa })_{\kappa \in -\bN} \} $ 
is a $ \hat{\mu} $-random walk indexed by $ -\bN $, i.e., 
\subitem (i) 
$ \hat{N}_{\kappa } $ is independent of $ \cF^{\hat{X},\hat{N}}_{\kappa -1} $ 
for all $ \kappa \in -\bN $; 
\subitem (ii) 
$ (\hat{N}_{\kappa })_{\kappa \in -\bN} $ is IID with common law $ \hat{\mu} $; 
{\subitem (iii) 
$ \hat{X}_{\kappa } = \hat{N}_{\kappa } \hat{X}_{\kappa -1} $ holds a.s. 
for all $ \kappa \in -\bN $. 
}

Moreover, it holds that 
\subitem (iv) 
$ \hat{X}_{\kappa } $ has uniform law on $ \hat{V} $ for all $ \kappa \in -\bN $. 
\end{Thm}

\Proof{
We have already shown (ii), (iii) and (iv). Let us prove (i). 
For this, it suffices to prove that, for any fixed $ \kappa \in -\bN $, 
\begin{align}
\text{$ \hat{N}_0,\hat{N}_{-1},\ldots,\hat{N}_{\kappa } $ and $ \hat{X}_{\kappa -1} $ 
are independent.} 
\label{}
\end{align}

Let $ B_{\kappa } \in \sigma(\hat{N}_0,\hat{N}_{-1},\ldots,\hat{N}_{\kappa }) $ 
and $ i \in \hat{V} $. 
Let $ l \in -\bN $. 
Since $ B_{\kappa } \cap \{ T_{\kappa -1} = l \} \in \sigma(N_0,\ldots,N_{l+1}) $, we have 
\begin{align}
& P \rbra{ B_{\kappa } , \ \hat{X}_{\kappa -1} = i , \ T_{\kappa -1} = l } 
\label{} \\
=& P \rbra{ B_{\kappa } , \ X_l \in \hat{V}_i , \ T_{\kappa -1} = l } 
\label{} \\
=& P \rbra{ B_{\kappa } , \ T_{\kappa -1} = l } P \rbra{ X_l \in \hat{V}_i } 
& \text{(by independence)} 
\label{} \\
=& \frac{1}{\hat{m}} P \rbra{ B_{\kappa } , \ T_{\kappa -1} = l } 
& \text{(by Theorem \ref{thm: main3})}. 
\label{}
\end{align}
This proves that 
$ \hat{X}_{\kappa -1} $ is independent of 
$ \sigma(\hat{N}_0,\hat{N}_{-1},\ldots,\hat{N}_{\kappa }) $. 

Let $ B_{\kappa +1} \in \sigma(\hat{N}_0,\hat{N}_{-1},\ldots,\hat{N}_{\kappa +1}) $ 
and $ \pi \in \fS(\hat{V}) $. 
Let $ l,l' \in -\bN $ with $ l > l' $. 
Since $ B_{\kappa +1} \cap \{ T_{\kappa } = l \} \in \sigma(N_0,\ldots,N_{l+1}) $, 
we obtain 
\begin{align}
& P \rbra{ B_{\kappa +1} , \ \hat{N}_{\kappa } = \pi , \ 
T_{\kappa } = l , \ T_{\kappa -1} = l' } 
\label{eq: prf line1} \\
= & P \rbra{ B_{\kappa +1} , \ T_{\kappa } = l , \ 
M \sbra{ N_l N_{l-1} \cdots N_{l'+1} } = \pi , \ 
T_{-1} \circ \theta_l = l' \Big. } 
\label{} \\
= & P \rbra{ B_{\kappa +1} , \ T_{\kappa } = l } 
P \rbra{ M \sbra{ N_l N_{l-1} \cdots N_{l'+1} } = \pi , \ 
T_{-1} \circ \theta_l = l' \Big. } 
\label{} \\
= & P \rbra{ B_{\kappa +1} , \ T_{\kappa } = l } 
P \rbra{ M \sbra{ N_0 N_{-1} \cdots N_{l'-l+1} } = \pi , \ 
T_{-1} = l'-l \Big. } 
\label{} \\
= & P \rbra{ B_{\kappa +1} , \ T_{\kappa } = l } 
P \rbra{ \hat{N}_0 = \pi , \ T_{-1} = l'-l \Big. } , 
\label{eq: prf line5}
\end{align}
where we write 
\begin{align}
T_{-1} \circ \theta_l = 
\max \{ k \le -p : N_{l+k+p} = \sigma_p , 
\ldots, N_{l+k+1} = \sigma_1 \} . 
\label{}
\end{align}
Summing up \eqref{eq: prf line1}-\eqref{eq: prf line5} by $ l' \in -\bN $, we have 
\begin{align}
P \rbra{ B_{\kappa +1} , \ \hat{N}_{\kappa } = \pi , \ 
T_{\kappa } = l } 
= P \rbra{ B_{\kappa +1} , \ T_{\kappa } = l } 
P \rbra{ \hat{N}_0 = \pi } . 
\label{}
\end{align}
This shows that $ \hat{N}_{\kappa } $ is independent of 
$ \sigma(\hat{N}_0,\ldots,\hat{N}_{\kappa +1}) $. 
Therefore, we conclude that 
$ \hat{N}_0,\hat{N}_{-1},\ldots,\hat{N}_{\kappa +1} $ and $ \hat{N}_{\kappa } $ 
are independent, 
which completes the proof. 
}

For $ k \in -\bN $, we write 
\begin{align}
K(k) = \max \{ \kappa \in -\bN : k-p \ge T_{\kappa } \} . 
\label{}
\end{align}
Note that 
\begin{align}
L(k) := T_{K(k)} 
= \max \{ l \le k-p : 
N_{l+p} = \sigma_p, \ldots, N_{l+2} = \sigma_2 , N_{l+1} = \sigma_1 \} 
\label{}
\end{align}
and that 
\begin{align}
\{ L(k)=l \} \in \sigma (N_k,N_{k-1},\ldots,N_{l+1}) 
, \quad l \le k-p . 
\label{}
\end{align}
The following theorem proves Claim \eqref{eq: main2}. 

\begin{Thm} \label{thm: main+}
Suppose that \eqref{eq: ass} holds. 
Then, for any $ k \in -\bN $, it holds that 
\begin{align}
X_k \in \cF^N_k \vee \sigma(\hat{X}_{K(k)}) 
\quad \text{a.s.} 
\label{eq: Xkin}
\end{align}
and that 
\begin{align}
\text{$ \hat{X}_{K(k)} $ is independent of $ \cF^N_0 $ 
and has uniform law on $ \hat{V} $.} 
\label{eq: hatXKk}
\end{align}
Consequently, 
if the road coloring is non-synchronizing, i.e., $ \hat{m} \ge 2 $, 
the $ \mu $-random walk $ \{ X,N \} $ is non-strong. 
\end{Thm}

\Proof{
By definitions of $ L(k) $ and $ K(k) $, we have 
\begin{align}
X_k 
=& N_k N_{k-1} \cdots N_{L(k)+p+1} N_{L(k)+p} \cdots N_{L(k)+1} X_{L(k)} 
\label{} \\
=& N_k N_{k-1} \cdots N_{L(k)+p+1} \sigma^{(n(1))} \cdots \sigma^{(n(p))} X_{L(k)} 
\label{} \\
=& N_k N_{k-1} \cdots N_{L(k)+p+1} x_i 
\label{}
\end{align}
with $ i= \hat{X}_{K(k)} $. 
Since $ N_k N_{k-1} \cdots N_{L(k)+p+1} \in \cF^N_k $, 
we obtain \eqref{eq: Xkin}. 

Let $ k' < k $ and let $ B \in \sigma(N_j:j \ge k'+1) $. 
Let $ l \le k-p $ and $ l' \le \min \{ k'-p,l \} $. 
Note that we have 
\begin{align}
& P \rbra{ X_{L(k)} \in \hat{V}_i, \ L(k)=l , \ L(k')=l' , \ B } 
\label{} \\
=& P \rbra{ X_l \in \hat{V}_i, \ L(k)=l , \ L(k')=l' , \ B } 
\label{} \\
=& P \rbra{ N_l \cdots N_{l'+p+1} \sigma_+ X_{l'} \in \hat{V}_i, \ L(k)=l , \ L(k')=l' , \ B } 
\label{} \\
=& \sum_{j=1}^{\hat{m}} P \rbra{ X_{l'} \in \hat{V}_j, \ M[N_l \cdots N_{l'+p+1}](j)=i , \ 
L(k)=l , \ L(k')=l' , \ B } . 
\label{eq: prf ind}
\end{align}
Since $ X_{l'} $ is independent of $ \cF^N_{0,l'} := \sigma(N_0,\ldots,N_{l'+1}) $, 
we see that 
\begin{align}
P \rbra{ \left. X_{l'} \in \hat{V}_j \right| \cF^N_{0,l'} } 
= P \rbra{ X_{l'} \in \hat{V}_j } 
= \frac{1}{\hat{m}} . 
\label{}
\end{align}
Since $ \{ L(k)=l, \ L(k')=l' \} \in \cF^N_{0,l'} $, we obtain 
\begin{align}
\text{\eqref{eq: prf ind}} 
=& \frac{1}{\hat{m}} \sum_{j=1}^{\hat{m}} P \rbra{ M[N_l \cdots N_{l'+p+1}](j)=i , \ 
L(k)=l , \ L(k')=l' , \ B \Big. } 
\label{} \\
=& \frac{1}{\hat{m}} P \rbra{ L(k)=l , \ L(k')=l' , \ B \Big. } . 
\label{}
\end{align}
Summing up in $ l \le k-p $ and $ l' \le k'-p $, we have 
\begin{align}
P \rbra{ \hat{X}_{K(k)} = i , \ B } 
= P \rbra{ X_{L(k)} \in \hat{V}_i , \ B } 
= \frac{1}{\hat{m}} P(B) . 
\label{}
\end{align}
This proves \eqref{eq: hatXKk}. 
}

\section{Periodic case} \label{sec: periodic}

If a directed graph $ (V,A) $ is strongly-connected, 
then it is easy to see that 
the greatest common divisor among $ \{ n \ge 1 : A^n(x,x) \ge 1 \} $ 
does not depend on $ x \in V $, 
so that it is called the {\em period} of $ (V,A) $. 
If the period of $ (V,A) $ is greater than one, then $ (V,A) $ is called {\em periodic}. 

Let us study periodic case. We shall utilize the following theorem. 

\begin{Thm}[Perron--Frobenius] \label{thm: PF2}
Let $ \mu $ be a probability law on $ \Sigma $. 
Suppose that 
the directed graph induced by $ \mu $ is strongly-connected and has period $ d \ge 2 $. 
Then there exist a partition $ \{ V^{(1)},\ldots,V^{(d)} \} $ of $ V $ 
and a family $ \{ \lambda^{(1)},\ldots,\lambda^{(d)} \} $ of probability laws on $ V $ such that 
the following assertions hold: 
\subitem {\rm (i)} 
for each $ x \in V^{(i)} $, it holds that 
\begin{align}
\mu(\sigma \in \Sigma : y=\sigma x) 
\begin{cases}
> 0 & \text{if $ y \in V^{(i+1)} $}, \\
= 0 & \text{otherwise}, 
\end{cases}
\label{}
\end{align}
where $ V^{(d+1)} = V^{(1)} $. 
\subitem {\rm (ii)} 
for each $ i=1,\ldots,d $, the support of $ \lambda^{(i)} $ is $ V^{(i)} $; 
\subitem {\rm (iii)} 
for each $ i=1,\ldots,d $ and each $ x,y \in V^{(i)} $, it holds that 
\begin{align}
\mu^{*nd}(\sigma \in \Sigma : y=\sigma x) 
\tend{}{n \to \infty } 
\lambda^{(i)}(y) . 
\label{}
\end{align}
\subitem {\rm (iv)} 
for each $ i=1,\ldots,d $, it holds that $ \mu^{*d} * \lambda^{(i)} = \lambda^{(i)} $. 
\end{Thm}

For the proof of Theorem \ref{thm: PF2}, 
see, e.g., \cite{MR2209438}. 
Each $ V^{(i)} $ will be called a {\em cyclic part}. 

In the sequel, 
let $ \mu $, $ \{ V^{(1)},\ldots,V^{(d)} \} $ and $ \{ \lambda^{(1)},\ldots,\lambda^{(d)} \} $ 
as in Theorem \ref{thm: PF2}. 
The following theorem chracterizes the class of all $ \mu $-random walks. 

\begin{Thm} \label{thm: periodic1}
The following assertions hold: 
\subitem {\rm (i)}
For each $ i=1,\ldots,d $, there exists a $ \mu $-random walk $ \{ X^{(i)},N^{(i)} \} $ in $ V $ 
such that $ X^{(i)}_0 \in V^{(i)} $ a.s. 
Such a $ \mu $-random walk is unique up to identity in law 
and the tail $ \sigma $-field $ \cF^{X^{(i)},N^{(i)}}_{-\infty } $ is trivial. 
\subitem {\rm (ii)}
Let $ \{ X,N \} $ be an arbitrary $ \mu $-random walk in $ V $. 
Then, for each $ i=1,\ldots,d $, the following hold: 
\subsubitem {\rm (ii-1)} 
the event $ A^{(i)} = \{ X_0 \in V^{(i)} \} $ belongs 
to the tail $ \sigma $-field $ \cF^{X,N}_{-\infty } $; 
\subsubitem {\rm (ii-2)} 
if $ P(A^{(i)})>0 $, then $ \{ X,N \} $ under $ P(\cdot \mid A^{(i)}) $ 
is identical in law to $ \{ X^{(i)},N^{(i)} \} $; 
in other words, it holds that 
\begin{align}
P((X,N) \in \cdot) 
= \sum_{i=1}^d P \rbra{ \rbra{ X^{(i)},N^{(i)} } \in \cdot } P(A^{(i)}) . 
\label{}
\end{align}
\end{Thm}

\Proof{
(i) 
For $ k=nd+r $ with $ n \in \bZ $ and $ r=1,2,\ldots,d $, 
we define $ \lambda^{(i)}_k = \mu^{*r} * \lambda^{(i)} $. 
By (iv) of Theorem \ref{thm: PF2}, we may apply Lemma \ref{lem: conv eq} 
to construct a $ \mu $-random walk $ \{ X^{(i)},N^{(i)} \} $ in $ V $ 
such that $ X^{(i)}_k $ has law $ \lambda^{(i)}_k $ for all $ k \in \bZ $. 
By (ii) of Theorem \ref{thm: PF2}, 
we have $ X^{(i)}_0 \in V^{(i)} $ a.s. 

Conversely, let $ \{ X^{(i)},N^{(i)} \} $ be a $ \mu $-random walk in $ V $ 
such that $ N^{(i)} $ has common law $ \mu $ and $ X^{(i)}_0 \in V^{(i)} $ a.s. 
For $ k \in \bZ $, let $ \lambda_k $ denote the law of $ X^{(i)}_k $. 
By (i) of Theorem \ref{thm: PF2}, we see that 
$ \lambda_{-nd} $ should have its support contained in $ V^{(i)} $ for all $ n \in \bN $. 
By (iii) of Theorem \ref{thm: PF2}, we see that 
\begin{align}
\lambda_0 = \mu^{*nd} * \lambda_{-nd} \tend{}{n \to \infty } \lambda^{(i)}. 
\label{}
\end{align}
This proves uniqueness by Lemma \ref{lem: conv eq}. 
The tail triviality can be proved in the same way as 
in the proof of Claim (v) of Theorem \ref{thm: main1}. 

(ii) 
Claim (ii-1) is obvious by (ii) of Theorem \ref{thm: PF2} 
and by the fact that 
$ X_k \in V^{(i)} $ if and only if $ X_{k-d} \in V^{(i)} $. 
Claim (ii-2) is obvious by the uniqueness in Claim (i). 
The proof is now complete. 
}

For $ i=1,\ldots,d $, 
we write $ \Sigma^{(i)} $ for the set of all mappings from $ V^{(i)} $ to itself, 
and define 
\begin{align}
\mu^{(i)} = \frac{1}{\mu^{*d}(\Sigma^{(i)})} \mu^{*d} |_{\Sigma^{(i)}} . 
\label{}
\end{align}
Then, by (i) of Theorem \ref{thm: PF2}, we see that 
$ \mu^{(i)} $ is a probability law on $ \Sigma^{(i)} $. 
The following theorem gives a necessary and sufficient condition 
for strongness of $ \mu $-random walks. 

\begin{Thm} \label{thm: periodic2}
Let $ \{ X,N \} $ be a $ \mu $-random walk in $ V $. 
\subitem {\rm (A)} 
If the tail $ \sigma $-field $ \cF^{X,N}_{-\infty } $ is trivial, 
then the $ \mu $-random walk $ \{ X,N \} $ is non-strong. 
\subitem {\rm (B)} 
Suppose that the tail $ \sigma $-field $ \cF^{X,N}_{-\infty } $ is non-trivial. 
Take $ i=1,\ldots,d $ such that $ \{ X,N \} \dist \{ X^{(i)},N^{(i)} \} $, 
which is possible by Theorem \ref{thm: periodic1}. 
Then following three assertions are equivalent: 
\subsubitem {\rm (i)} 
$ \Supp(\mu^{(i)}) $ is synchronizing. 
\subsubitem {\rm (ii)} 
The limit $ \displaystyle \lim_{l \to -\infty } N_k N_{k-1} \cdots N_{l+1} $ 
exists a.s. for all $ k \in \bZ $. 
\subsubitem {\rm (iii)} 
The $ \mu $-random walk $ \{ X,N \} $ is strong.  
\end{Thm}

The proof of Theorem \ref{thm: periodic2} is an immediate consequence 
of Theorem \ref{thm: main intro}, 
so that we omit it.

{\bf Acknowledgements:} 
The author thanks Kenji Yasutomi 
for his suggestions about Lemma \ref{thm: F-clique} and the example in Section \ref{sec: ex3}. 
The author also thanks Professors Yusuke Higuchi and Tomoyuki Shirai 
for fruitful discussions.

\def\cprime{$'$} \def\cprime{$'$}

\end{document}